\documentclass[letterpaper, 10pt, onecolumn,draftcls]{IEEEtran} 
\pdfoutput=1
\usepackage[margin=1in]{geometry}
\usepackage{amsmath,amsfonts,amssymb,mathtools}
\usepackage{bbm}
\usepackage{color}
\usepackage{graphicx}
\usepackage{enumerate}
\usepackage{mathtools}
\mathtoolsset{showonlyrefs=true}
\newcommand{\R}{\mathbb{R}}

\newcommand{\bone}{\mathbbm{1}}
\newtheorem{definition}{Definition}

\newtheorem{lemma}{Lemma}

\newtheorem{theorem}{Theorem}
\newtheorem{assumption}{Assumption}

\newcommand{\diag}[2]{\textnormal{diag}(#1, \ldots, #2)}

\author{Matthew T. Hale$^\star$, Sebastian F. Ruf$^\dagger$, Talha Manzoor$^\ddagger$, Abubakr Muhammad$^\S$
\thanks{$^\star$Department of Mechanical
and Aerospace Engineering, University of Florida, Gainesville, FL, USA. Email: \texttt{matthewhale@ufl.edu}.}
\thanks{$^\dagger$Center for Complex Network Research and Department of Psychology,
Northeastern University, Boston, MA, USA. Email: \texttt{sebastianfruf@gmail.com}.}
\thanks{$^\ddagger$Department of Electrical Engineering, Namal College, Mianwali, Pakistan. Email: \texttt{talha@namal.edu.pk}.}
\thanks{$^\S$Department of Electrical Engineering,
Lahore University of Management Science, Lahore, Pakistan. Email: \texttt{abubakr@lums.edu.pk}.}
}

\begin{document}
\title{Stability and Sustainability of Resource Consumption Networks}
\maketitle

\begin{abstract}
In this paper, we examine both stability and sustainability of 
a network-based model of natural resource
consumption. Stability is studied from a dynamical systems perspective, though 
we argue that sustainability is a fundamentally different notion from stability in social-
ecological systems. Accordingly, we also
present a criterion for sustainability that is guided by the existing literature on 
sustainable development. 
Assuming a generic social network of consuming agents' interactions, 
we derive sufficient conditions for both 
the stability and sustainability of the model as constraints on the network structure itself. 
We complement these 
analytical results with numerical simulations and discuss the implications of our findings for 
policy-making for sustainable resource governance.
\end{abstract}
\section{Introduction}
\label{sec:intro}

Sustainability of social-ecological systems has been a subject of considerable interest for some 
time~\cite{ostrom2015governing}. However, despite its importance, 
sustainability of this kind is still far from being adequately defined in a rigorous context. The confusion stems directly from the definition of sustainable development itself, given by the Brundtland Commision of the World Commission on Environment and Development (WCED)~\cite{brundtland1987our}, which leaves much room for interpretation. This has led to a string of research spanning across multiple disciplines in search of a formal definition of the concept~\cite{cabezas2002towards, chichilnisky1995green, martinet2009defining}. 

From a dynamical systems perspective, sustainability has often been linked to developments 
in~\cite{ludwig1997sustainability, patten1997logical}, and sometimes even used interchangeably 
with stability in social-ecological systems~\cite{kinzig2014consumption}. 
We argue that while stability is indeed relevant in this setting, it is a purely dynamical systems property that does not appropriately account for the ecological nature of the system. Thus a separate formulation of sustainability is needed. This need is underscored by the rising interest of the controls community not only in social-ecological systems, but also in other systems lying at the interface of technological, ecological, and social sciences. Such systems are collectively termed
``cyber-physical social systems''~\cite{wang2010emergence}.

In this paper we study the stability and sustainability of a social-ecological system whose model was first presented and studied in \cite{manzoor2016game}. The model is for a natural resource being harvested by a fixed number of consuming agents. While deciding their consumption, the agents not only take information about the resource into account, but also incorporate information about the consumption of other neighboring agents. The social network of the population thus plays a critical role in determining overall system behavior, and 
it must also affect sustainability in any sense that is considered. Indeed, it has been 
found~\cite{ostrom2015governing} that the ability of a population to sustainably manage its resources at the community level is highly dependent on the underlying social network of that community. 

Accordingly, the insights provided by the network structure 
have resulted in an increase in the application of network analysis tools to natural resource 
systems~\cite{prell2011social, videras2013social}. For instance, social structure has been observed to be strongly correlated with pro-environmental behavior in social-ecological settings \cite{videras2012influence}. Information about the structure of the network has also been exploited in the past to identify elements critical to sustainable governance of a resource (for some representative studies, see~\cite{prell2009stakeholder,crowe2007search,ramirez2009impact}). Previously, our selected model of resource consumption has been studied only under restrictive assumptions on the network topology~\cite{manzoor2018learning}. 

In this paper, we maintain a generic network structure  to fully understand the effects of the social 
network on sustainability of its resource. A complete characterization of sustainable communities 
 remains an open problem for the scientific community. This has been attributed in part to the highly 
 inter-disciplinary nature of the field and a lack of integrative studies to unify the research that 
 remains scattered among different strands of work~\cite{rockenbauch2017social}. In this paper we 
 devise a criterion for sustainability that is subsequently applied to our network model of resource 
 consumption. This results in structural conditions on the network topology which, we believe, 
 contributes a step towards uncovering the structural characteristics of sustainable societies.

In what follows, we develop a sustainability criterion that draws from the literature on sustainable development. This criterion is based on a rigorous definition of sustainability that we introduce, 
and we show how the enforcement of system sustainability is translated into conditions on the network topology and its associated parameters. 
We apply this criterion to a network model of resource-consuming agents and compare it with  sufficient conditions for stability. To that end, we also present a
novel proof of stability for the full $n$-agent model which was previously known to be stable only for 2 agents \cite{manzoor2018learning}. The contributions of this paper thus consist of
the stability proof for the general~$n$-agent model and the sustainability criterion we
develop, together with the resulting analyses and implications
for resource consumption networks that we discuss. 


The remainder of the paper is organized as follows. Section~\ref{sec:model} first provides the model 
of interest. 
Then, Section~\ref{sec:stability} presents sufficient conditions for stability of the model
and a Lyapunov-based stability proof. 
Section~\ref{sec:sust} next provides a commentary on the sustainability literature originating from different disciplines. Following guidelines from this exposition, we also present our mathematical definition of sustainability. In Section~\ref{sec:proof}, we obtain
sufficient conditions for system sustainability
in the form of constraints on the structure of its
underlying social network. Section~\ref{sec:sim} presents simulations of both stable and 
sustainable networks along with a discussion on the implications for policy making. 
We conclude in Section~\ref{sec:conc}.

\section{Background and Network Model} \label{sec:model}
This section describes the system model of interest and its coupled resource and consumption dynamics. The model was first presented in \cite{manzoor2016game} and subsequently studied in \cite{manzoor2018learning,manzoor2017structural,ruf2018stability}. For full details on the environmental and social-psychological foundations of the model, we refer the reader to the sources mentioned above. Here we present a non-dimensionalized version of the model (also introduced in \cite{manzoor2016game}) that not only has a reduced parameter space, but also allows more straightforward interpretations of the results relative to the original model. 

The setting is that of a single natural resource whose stock renews according to the standard model of logistic growth \cite{perman2003natural}. 
The resource has an associated carrying capacity and intrinsic growth rate that affect
the evolution of its value over time. 
 Let $x(t)>0$ represent the resource quantity at time $t$ relative to its carrying capacity. Thus $x(t)=1$ implies that at time $t$ the resource stock is at its environmental carrying capacity. The resource is harvested by a consuming population consisting of~$n$ agents. Each agent harvests the resource by exerting some effort, and,
in this model, 
 this effort constitutes all the different dimensions of harvesting activity. In fishing for instance, it can capture the number of boats deployed, their efficiency, the number of days fishing is undertaken, 
etc.~\cite{perman2003natural} 

Let $y_i(t) \in \mathbb{R}$ represent the consumption effort 
of agent~$i$ relative to the rate of growth of the resource; 
a consumption effort of $y_i(t)=1$ implies a 
consumption 
rate that is equal to the intrinsic growth rate of the resource. The resource dynamics are then given by 
\begin{align}
	\dot{x}(t) =  (1-x(t))x(t) - x(t)\sum_{i=1}^{n}  y_i(t),
\end{align}
where $i \in \{1,\dots ,n\}$ indexes the agents in the network.  
Note that the model permits a negative consumption rate. From the dynamics of the resource it is evident that this may result in $x(t)$ taking on values greater than unity, which corresponds to the stock level crossing the natural carrying capacity of the environment. 
In \cite{manzoor2016game} we discuss the physical interpretations of this phenomenon in detail. 
While a positive
consumption rate represents an effort to consume the resource (thereby decreasing its stock), 
a negative rate corresponds to an action 
to increase the stock of the resource. 
This includes actions that directly increase the stock (such as planting trees or breeding fish) as well as more indirect actions (such as restoration of soil fertility or a favorable distribution of gear to local fishermen). These actions may also be directed towards artificially increasing the carrying capacity, e.g., de-silting water canals, increasing the land available for forest growth, shifting to more intensive farming techniques, etc.

Next we present the dynamics of the individual consumption efforts $y_i(t)$. According to  Festinger's theory of social comparison processes~\cite{festinger1954theory}, while making decisions, humans incorporate both objective information and social information. In the context of resource consumption \cite{mosler2003integrating}, the objective information is interpreted as information on the state of the resource, i.e., whether it is abundant or scarce. We model this through the scarcity threshold $\rho_i > 0$ associated with agent~$i$. The state of the resource is then captured by the ecological factor $(x(t) - \rho_i(t))$. A negative ecological factor indicates that
agent~$i$ perceives the resource as scarce, whereas a positive factor indicates 
that agent~$i$ perceives the resource to be abundant. The social information is interpreted as the difference between agent~$i$'s consumption level and the consumption
levels of other neighboring agents. This is represented by the social factor $\sum_{j=1}^{n} \omega_{ij}\left( y_j(t) - y_i(t) \right)$, where $\omega_{ij} > 0$ is the directed tie-strength between
agents~$i$ and $j$ and represents the influence that agent~$j$ has on agent~$i$'s consumption. Furthermore $\omega_{ii} = 0$ and $\sum_{j=1}^{N} \omega_{ij} = 1$ for all~$i$. 

The rate of change of consumption effort is then represented as a weighted sum of the ecological and social factors. These are weighed by the ecological and social 
weights ${\alpha_i \in [0,1]}$ and ${\nu_i \in [0,1]}$, respectively (these terms are also referred to as the social and ecological relevances of agent~$i$ below). Since the weighing of both factors in the final decision-making is observed to be one-dimensional~\cite{festinger1954theory}, the weights  sum to one, i.e., $\alpha_i + \nu_i = 1$ for all $i$. Thus the consumption dynamics are given by 
\begin{equation}
	\dot{y}_i(t) = b_i \Big( \alpha_i(x(t) -\rho_i)- \nu_i \sum_{j=1}^{n} \omega_{ij}\left( y_i(t) - y_j(t) \right) \Big),
\end{equation}
where $b_i > 0$ is called the sensitivity of 
agent~$i$ and represents the readiness of agent~$i$ to change her consumption.

Based on the above, the coupled dynamics of the resource stock and consumption effort are
\begin{align}\label{eq:ses}
	\dot{x}(t) &=  (1-x(t))x(t) - x(t)\sum_{i=1}^{n}  y_i(t),\\
    \dot{y}_i(t) &= b_i \Big( \alpha_i(x(t) -\rho_i)- \nu_i \sum_{j=1}^{n} \omega_{ij}\left( y_i(t) - y_j(t) \right) \Big),
\end{align}
where~$i \in \{1, \ldots, n\}$. 
In Equation~\eqref{eq:ses}, $t$ is non-dimensional time which is normalized with respect to the resource growth rate. In~\cite{manzoor2016game}, we have discussed how this normalization 
reduces the dimension of the parameter space relative to the original model, 
along with interpretations of the various variables and their domains, and we refer the reader to that
reference for an extended discussion of the subject.  
\section{Stability of Networked Resource Consumption} \label{sec:stability}
This section proves that the network consumption model
in Section~\ref{sec:model} is asymptotically stable. We first rewrite
the dynamics as an equivalent aggregate-level resource
and consumption system and then transform this system
into one whose equilibrium is at the origin. Then
we provide a Lyapunov-based stability proof. Below, we use 
$\diag{a_1}{a_n}$ to denote the diagonal
matrix with the scalars $a_1$ through $a_n$ on its main diagonal. 

\subsection{Aggregate Network Dynamics}
We define the states~$\gamma = \ln x$ and~$y = (y_1, \ldots, y_n)^T$ 
along with the matrices
$A = \diag{\alpha_1}{\alpha_n}$,
${B = \diag{b_1}{b_n}}$, and
$V = \diag{\nu_1}{\nu_n}$. We further define
the vector $\rho = (\rho_1, \ldots, \rho_n)^T$
and the matrix of weights
\begin{equation}
T = \left(\begin{array}{ccccc} 
1 & -\omega_{12} & -\omega_{13} & \cdots & -\omega_{1n} \\
-\omega_{21} & 1 & -\omega_{23} & \cdots & -\omega_{2n} \\
\vdots & \vdots & \vdots & \ddots & \vdots \\
-\omega_{n1} & -\omega_{n2} & -\omega_{n3} & \cdots & 1
\end{array}\right).
\end{equation}
We further define $\theta_i = \frac{\nu_i}{\alpha_i}$, which we will
use below.

Differentiating $\gamma$ and $y$ with respect to time gives
\begin{align}
\dot{\gamma} &= \frac{\dot{x}}{x} = 1 - x - \bone^Ty \\
\dot{y} &= BA(e^{\gamma}\bone - \rho) - BVTy. 
\end{align}
In this form, we impose the following assumptions.

\begin{assumption} \label{as:inv}
The matrix $(A\bone\bone^T + VT)^{-1}$ exists and
\begin{equation}
\bone^T\big(A\bone\bone^T + VT\big)^{-1}A(\bone - \rho)  < 1. \tag*{$\triangle$}
\end{equation}
\end{assumption}
 Below, Assumption~\ref{as:inv} will be used to insure that $\gamma_0$ and~$y_0$,
the respective equilibria of the~$\gamma$ and~$y$ dynamics, are well-defined. 
\begin{assumption} \label{as:conn}
The underlying graph connecting agents is strongly connected. \hfill $\triangle$
\end{assumption}
Assumption~\ref{as:conn} is a common assumption in multi-agent networks. Strong
connectivity implies, roughly, that there is ``enough'' interaction among agents
to drive the state of the system toward its equilibrium point, and this notion
will be made precise in Theorem~\ref{thm:stability}. 
\begin{assumption} \label{as:weights}
For all agent indices $i \in \{1, \ldots, n\}$,
\begin{equation}
\theta_i \geq \sum_{\substack{k = 1 \\ k \neq i}}^{n} \omega_{ki}\theta_k. \tag*{$\triangle$}
\end{equation} 
\end{assumption}
Because~$\theta_i = \frac{\nu_i}{\alpha_i}$, we can regard~$\theta_i$ as measure of how
social agent~$i$ is. Assumption~\ref{as:weights} then requires that agent~$i$'s outgoing
influence, measured upon agent~$k$ through~$\omega_{ki}$, be limited based on its
own sociability. To influence other highly social agents, agent~$i$ must itself
be more social. Conversely, if agent~$i$ is less social, then it will have a smaller
value of~$\theta_i$, and the weights of its outgoing edges must be smaller. 

Computing the equilibrium of these dynamics, we set~$\dot{\gamma} = 0$,
which gives~$\gamma_0 = \ln(1 - \bone^Ty_0)$. 
To find $y_0$, we have
$\dot{y} = 0 = B(A\bone e^{\gamma_0} - A\rho - VTy_0)$ 
where solving for $y_0$ gives
$y_0 = (A\bone\bone^T + VT)^{-1}A(\bone - \rho)$,
which is well-defined under Assumption~\ref{as:inv}. 

To shift this system's equilibrium to the origin,
we then define the states~$v = \gamma - \gamma_0$
and~$w = y - y_0$. 
Computing $\dot{v}$ and $\dot{w}$, we find 
\begin{subequations} \label{eq:vwdynamics}
\begin{align} 
\dot{v} &= -e^{\gamma_0}(e^v - 1) - \bone^Tw \label{eq:v} \\
\dot{w} &= e^{\gamma_0} BA\bone(e^v - 1) - BVTw.\label{eq:w} 
\end{align}
\end{subequations}

\subsection{Global Asymptotic Stability}
The following theorem shows that  the system model in Equation~\eqref{eq:vwdynamics} is 
globally asymptotically stable. Due to its technical nature, we relegate
the proof to the appendix. 

\begin{theorem} \label{thm:stability}
The system in Equation~\eqref{eq:vwdynamics} is globally asymptotically stable to the origin.
\end{theorem}
\emph{Proof:} See the Appendix. \hfill $\blacksquare$

While stability is of course critical to insuring acceptable asymptotic network behaviors, 
stability as a purely dynamical systems notion does
not account for the fact that this is a resource consumption network.
The value of the resource itself must be sustained over time, and thus
its dynamics and those of the agents' consumptions must be accounted for
not only in the limit, but also at each point in time as the system
evolves. Accordingly, we next assess the role of sustainability in this
system.

\section{Formalizing Sustainability}\label{sec:sust}
In the quest for attaining sustainable societies, one of the biggest challenges is formalizing a notion of sustainability. The  
Brundtland Commission of the WCED~\cite{brundtland1987our} defines sustainable development as ``development that meets the needs of the present without compromising the ability of future generations to meet their own needs''. Clearly, this definition leaves much room for interpretation, and a vast amount of literature is devoted to developing a definitive measure of sustainability. As a result, there exist varying notions of sustainability originating from different disciplines ranging from the ecological sciences \cite{cabezas2002towards} to economics \cite{pezzey2017economics}. These discrepancies 
create significant 
difficulty in characterizing sustainable behavior in the rigorous framework of dynamical systems.

It is quite natural to think of sustainable consumption in social-ecological systems as a state in which the level of consumption is allowed to be non-decreasing~\cite{valente2005sustainable}. On the other hand, 
it is broadly understood that a continuously increasing consumption is not realizable in a finite world with limited resources~\cite{meadows2012limits}. Sustainability would then imply an indefinitely maintainable level of consumption, which corresponds to attaining a positive equilibrium. It is not surprising therefore, that stability has been closely linked with sustainability by researchers working in the area (see for instance \cite{ludwig1997sustainability,patten1997logical,kinzig2014consumption}). 

However, while stability is roughly understood to be a necessary property for sustainability, it is not sufficient~\cite{patten1997logical}. 
In particular, from an economical viewpoint, Chichilnisky~\cite{chichilnisky1995green} argues that sustainability should not only rely on the behavior of the consumption path in the limit (as time approaches infinity), but also on the accumulated utility for all time. From an ecological perspective, an ecosystem can remain functional only if its critical parameters remain within a certain region over time~\cite{ben2012cybernetics}, which likewise
suggests the use of some non-asymptotic analysis.

Martinet formalizes this phenomenon in his definition for sustainability as bounds on certain sustainability 
indicators~\cite{martinet2009defining}. 
He then defines a
sustainability criterion as a generalized maximin problem in search of the bounds that are optimal 
with respect to certain preference functions. 
The boundedness of system response is captured by the concept of Lagrange stability  in dynamical systems theory~\cite{bhatia2002stability}. Lagrange stability has also been related in the ecological literature  with the ability of an ecosystem to resist disturbances~\cite{patten1997logical}. However Lagrange stability (along with other classical notions of stability \cite{leine2010historical}) is an absolute property of the system
rather than its trajectories. 
This poses a challenge in prioritizing different trajectories based on their level of sustainability. In other words, stability alone does not capture the nuances of sustainability. 

Beyond stability, 
another property of the system that is closely linked with sustainability is resilience. 
Holling~\cite{holling1973resilience} defines resilience as the ability of a system to maintain its integrity when subjected to disturbances. From a dynamical systems viewpoint, this property can be related to how rapidly or gradually a system moves towards an equilibrium. 
For social-ecological systems, a sudden change in consumption patterns may trigger a major regime change in the underlying habitat accommodating the natural resource, which in turn may cause unpredictable
and undesirable changes in the resource stock itself. 
Thus, a consumption path must not vary at more than a certain rate for it to be  sustainable. The notion of sustainability as permitting
a non-decreasing consumption path already implies a lower bound on the derivative of the
consumption trajectory~\cite{loucks1997quantifying}. 
However the above argument also implies an upper bound on the derivative. It is important to note that while a limit on the rate of change of the trajectories protects the system from loss of resilience, it also further reflects the physical limits associated with the rate of resource growth and harvesting.

Even across the different notions of sustainability present in the 
literature, 
e.g.,~\cite{valente2005sustainable,loucks1997quantifying,ulph2014sustainable,kharrazi2013quantifying},
the underlying philosophy remains the same: sustainability entails preserving essential system characteristics over a sufficiently long period of time. What distinguishes these notions then, is the different perspectives on what characteristics of the system must persist and for how long. Costanza and Patten~\cite{patten1997logical} assert that sustainability is highly dependent on the scale at which the system is perceived. 
In particular, the sustainability of a larger system does not imply the sustainability of all of its subsystems. For example, cells and bacteria die periodically in order for a larger organism to sustain its life,
or a small city may be consumed by construction of a dam in order to promote the sustainability of the region as a whole. The losses
of these smaller subsystems prolong the sustainability of the larger systems, though it is broadly
understood that this does not provide ``maintenance forever''
of the larger system~\cite{costanza1995defining}.
Sustainability thus cannot imply an infinite lifetime. Rather, it implies that the system achieves its expected lifespan, and this lifespan depends both on the temporal and spatial scales of the system of 
interest~\cite{costanza1995defining}.

In light of the above discussion, it is possible to identify the salient features that we 
can use to mathematically encode sustainability. These are listed in sequence below: 

\begin{itemize}
	\item First, sustainability should be evaluated over finite intervals. This reflects that sustainability of a system does not necessarily mean indefinite survival, but rather the attainment of its expected 
	lifespan~\cite{patten1997logical}. As discussed above, this lifespan depends on the natural temporal and spatial scales of a system. 
	\item Second, sustainability should provide criteria pointwise in time as opposed to only at the steady-state
	behavior of a system. Sustainability concerns itself not only with the generation in the time limit, but also with the current generation and all others in between~\cite{brundtland1987our}. Thus, rather than being asymptotic, a definition of sustainability should account for all states of a system over time \cite{chichilnisky1995green,valente2005sustainable}.
	\item Third, sustainability should be a structural property of a system~\cite{ben2012cybernetics}. From a systems perspective, the structure of the system determines the characteristics that are essential to attaining a particular state of being. In this context, sustainability may be thought of as the invariance of a particular quantity (mediated by the internal structure) that captures the required behavior. The model in Section~\ref{sec:model} is a networked system, and the ``structure'' of this system consists of both its topology and associated edge weights, which are encapsulated by the matrix~$T$. 
\end{itemize}
Following the guidelines listed above, we develop a mathematical formulation of sustainability as follows. First,
we consider sustainability for systems across the time horizon $[0, t_{max}]$ for some~$t_{max} > 0$. 
Second, we define the following four constants: 
\begin{enumerate}
	\item $v_{max}$, the maximum allowable value of $v$
	\item $v_{min}$, the minimum allowable value of $v$
	\item $d_{max}$, the maximum allowable value of $\dot{v}$
	\item $d_{min}$, the minimum allowable value of $\dot{v}$,
\end{enumerate}
where~$d$ is chosen to indicate ``derivative''. 

Together, the four constants listed above define a box in 
the $(v, \dot{v})$-plane, and sustainability of a system is equivalent to invariance of this box. We emphasize that the equilibrium of~$v$ need not be in this box, which further highlights the distinction between stability and sustainability. 
Sustainability can then be formally expressed as an invariance condition. 
We define the subset $S$ of the $(v,\dot{v})$-plane as
\begin{equation}
	S := \left\{\big(v, \dot{v}\big) \mid v \geq v_{min}, v \leq v_{max}, \dot{v} \geq d_{min}, \dot{v} \leq d_{max}\right\},
\end{equation}
and sustainability for a social-ecological system is then defined as follows.

\begin{definition} \label{def:sust}
	A quantity of resource $v(t)$ is sustainable over  $[0, t_{max}]$ 
	if $\big(v(t), \dot{v}(t)\big) \!\in\! S$ for all $t \!\in\! [0, t_{max}]$. \hfill $\triangle$
\end{definition}

Note that many of the above sustainability notions from the existing literature
pertain directly to consumption, while our
mathematical sustainability criterion is stated in terms of the resource. 
Despite this apparent difference, we show in Section~\ref{sec:proof} below that
Definition~\ref{def:sust} must account for agents' consumptions to be satisfied.
In particular, we show how Definition~\ref{def:sust} 
is linked with both the social structure and individual parameters of the consuming network by deriving conditions on $T, A$ and $B$ that are sufficient to satisfy 
Definition~\ref{def:sust}.

\section{Insuring Sustainability of Networked Resource Consumption}
\label{sec:proof}
The social structure of a resource consumption network is specified
both by the graph of agents' interactions and the weights
of edges in this graph. These values are captured
in the matrix~$T$ in Equation~\eqref{eq:vwdynamics}, and
this section will derive conditions on $T$ that insure that
$v$ remains at sustainable levels over some pre-specified
time horizon $[0, t_{max}]$. These results require not
only accounting for how $v$ evolves, but also how agents' 
consumption levels affect $v$. 

\subsection{Preliminaries}
Towards deriving the desired conditions on $T$, we
first give some preliminary results related to $v$
and $w$ that we will use below. We begin by stating
the Bellman-Gr\"{o}nwall inequality in the required form. 

\begin{lemma} \label{lem:bellman} \emph{Bellman-Gr\"{o}nwall Inequality; \cite{ames1997inequalities})} 
Let $K_1 \in \R$ and $K_2 > 0$. If, for all~$t \in [a, b]$, $\phi : [a, b] \to \R$, the bound
$\phi(t) \leq K_1 + K_2\int_{a}^{t} \phi(s) \, ds$ is satisfied,
then on the same interval
$\phi(t) \leq K_1e^{K_2 (t-a)}$. \hfill $\blacksquare$ 
\end{lemma}

Next, we enforce the following assumption.

\begin{assumption} \label{as:v}
We assume that
\begin{itemize}
\item $v_{min} < v_{max}$
\item $v_{max} > 0$
\item $v_{min} < v(0) < v_{max}$
\item $d_{min} < 0 < d_{max}$. \hfill $\lozenge$
\end{itemize}
\end{assumption}

This assumption is rather mild in what it imposes upon the system. Requiring 
that~$v_{max} > 0$ and~$v_{max} > v_{min}$ insures that a given
set of sustainability bounds is feasible. Assuming that~$v_{min} < v(0) < v_{max}$
merely assumes that a system starts within the bounds it must remain within.
And requiring that~$d_{min} < 0 < d_{max}$ means that the resource levels in a
system should be allowed to both increase and decrease, which is 
natural. 

Next, we define the constants
$\beta = \max_{i \in \{1, \ldots, n\}} b_i\nu_i$, 
$C_1 = \|w(0)\|_1 + t_{max}e^{\gamma_0}(e^{v_{max}} - 1)\sum_{i=1}^{n} b_i\alpha_i$, and 
$C_2 = \beta \|T\|_1t_{max}$,
and we will repeatedly encounter them below. 
Using Lemma~\ref{lem:bellman} and Assumption~\ref{as:v}, we have
the following characterization of $w$.
\begin{lemma} \label{lem:normw}
Suppose that $v(s) \leq v_{max}$ for all $s \in [0, t)$. 
Then~$\|w(t)\|_1 \leq C_1\exp(C_2)$.
\end{lemma}
\emph{Proof:} The fundamental theorem of calculus
gives
\begin{equation}
w(t) = w(0) + \int_{0}^{t} e^{\gamma_0}BA\bone(e^v - 1) - BVTw(\tau)d\tau.
\end{equation}
Taking the $1$-norm of both sides and applying the triangle inequality
then gives
\begin{align}
\|w(t)\|_1 &\leq \|w(0)\|_1 + \int_{0}^{t} \|e^{\gamma_0}BA\bone(e^{v(\tau)} - 1)\|_1 d\tau
+ \int_{0}^{t} \|BV\|_1\|T\|_1\|w(\tau)\|_1 d\tau \\
   &\leq \|w(0)\|_1 + t_{max}e^{\gamma_0}\|BA\bone(e^{v_{max}} - 1)\|_1 + \|BV\|_1\|T\|_1\int_{0}^{t} \|w(\tau)\|_1 d\tau \\
   &= \|w(0)\|_1 + t_{max}e^{\gamma_0}(e^{v_{max}} - 1)\sum_{i=1}^{n} b_i\alpha_i + \|BV\|_1\|T\|_1\int_{0}^{t} \|w(\tau)\|_1 d\tau,
\end{align}
where we have used that $v(\tau) \leq v_{max}$ and expanded the
$1$-norm of $BA\bone$. 

Using that the $1$-norm of a matrix is equal to its largest column sum, the
value of $\|BV\|_1$ is then equal to $\beta$ because $BV$ is diagonal.
The lemma follows by applying
Lemma~\ref{lem:bellman} with
$K_1 = C_1$ and 
$K_2 = C_2$. 
\hfill $\blacksquare$

We next derive a similar bound for $v$.

\begin{lemma} \label{lem:normv}
Suppose that $v(s) \leq v_{max}$ for all $s \in [0, t)$. Under the dynamics
in Equation~\eqref{eq:vwdynamics}, we find that 
\begin{equation}
v(t) \geq v(0) - t_{max}\big[e^{\gamma_0}(e^{v_{max}} - 1) - C_1\exp(C_2)\big].
\end{equation}
\end{lemma}
\emph{Proof:} From the fundamental theorem of calculus, 
\begin{equation}
v(t) = v(0) - \int_{0}^{t} e^{\gamma_0}(e^{v(s)} - 1)ds - \int_{0}^{t} \bone^Tw(s) ds. 
\end{equation}
Noting that~$e^{v(\tau)} - 1 \leq e^{v_{max}} - 1$,
we find 
\begin{equation}
v(t) \geq v(0) - te^{\gamma_0}(e^{v_{max}} - 1) 
- \int_{0}^{t} \bone^Tw(\tau) d\tau. 
\end{equation}

For the final term above we find that
\begin{equation}
\int_{0}^{t} \bone^Tw(\tau)d\tau 
\leq \int_{0}^{t} \|w(\tau)\|_1 d\tau. 
\end{equation}
Multiplying by $-1$ to reverse the inequality, we find
\begin{equation}
v(t) \geq v(0) - te^{\gamma_0}(e^{v_{max}} - 1) 
- \int_{0}^{t} \|w(\tau)\|_1 d\tau. 
\end{equation}
We complete the proof by applying Lemma~\ref{lem:normw}. \hfill $\blacksquare$

Having established
these basic preliminary lemmas,
we next derive bounds on $T$ to enforce
the sustainability bounds on $v$ and $\dot{v}$. 

\subsection{Enforcing $v(t) \leq v_{max}$}
To enforce $v(t) \leq v_{max}$, we enforce 
$\dot{v}\big|_{v = v_{max}} \leq 0$,
and the following lemma gives a sufficient condition for doing so.

\begin{lemma} \label{lem:vmax}
Under Assumption~\ref{as:v}, 
if $T$ satisfies
\begin{equation}
\|T\|_1 \leq \frac{1}{\beta t_{max}}\ln\left(\frac{e^{\gamma_0}(e^{v_{max}} - 1)}{C_1}\right)
\end{equation}
then $v(t) \leq v_{max}$ for all $t \in [0, t_{max}]$. 
\end{lemma}
\emph{Proof:} 
Using Equation~\eqref{eq:v}, 
enforcing~$\dot{v}\Big|_{v = v_{max}} \leq 0$
is equivalent to requiring
$-e^{\gamma_0}(e^{v_{max}} - 1) - \bone^Tw \leq 0$,
which we rearrange to find
$-\bone^Tw \leq e^{\gamma_0}(e^{v_{max}} - 1)$.

We note that
$-\bone^Tw \leq \|w\|_1$,
and we will thus enforce the sufficient condition
\begin{equation} \label{eq:l4sufficient}
\|w(t)\|_1 \leq e^{\gamma_0}(e^{v_{max}} - 1)
\end{equation}
for all $t$. 
Using Assumption~\ref{as:v},~$v(0) < v_{max}$. Suppose that~$v$ approaches~$v_{max}$
at time~$t$, i.e.,~$v(t) = v_{max}$. Then the condition~$v(s) \leq v_{max}$
for all~$s \in [0, t)$ is satisfied and Lemma~\ref{lem:normw} can be used at time~$t$.
Then a sufficient condition for enforcing Equation~\eqref{eq:l4sufficient}
is given by
\begin{equation}
C_1\exp(\beta\|T\|_1 t_{max}) \leq e^{\gamma_0} (e^{v_{max}} - 1).
\end{equation}
Solving for~$\|T\|_1$ then gives a bound that insures~$\dot{v}(t)\big|_{v = v_{max}} \leq 0$.
If~$v$ later approaches~$v_{max}$ at some time~$\tau$, then
repeating the above argument at~$\tau$ insures that~$\dot{v}(\tau)\big|_{v = v_{max}} \leq 0$ as well.
Of course, if~$v$ never approaches~$v_{max}$ then~$v(t) \leq v_{max}$ trivially holds, and thus
it is true for all~$t \in [0, t_{max}]$ in all cases. 
\hfill $\blacksquare$

\subsection{Enforcing $v(t) \geq v_{min}$}
We next derive a sufficient condition on $T$ that implies
that the bound $v(t) \geq v_{min}$ is satisfied for all $t \in [0, t_{max}]$. 

\begin{lemma} \label{lem:vmin}
If $\|T\|_1$ satisfies the bound in Lemma~\ref{lem:vmax} and 
\begin{equation}
\|T\|_1 \leq \frac{1}{\beta t_{max}}\ln\left(\frac{v(0) - t_{max}e^{\gamma_0}(e^{v_{max}} - 1) - v_{min}}{t_{max}C_1}\right),
\end{equation}
then $v(t) \geq v_{min}$ for all $t \in [0, t_{max}]$. 
\end{lemma}
\emph{Proof:} 
The bound from Lemma~\ref{lem:vmax} insures that~$v(t) \leq v_{max}$ for 
all~$t \in [0, t_{max}]$ and thus Lemma~\ref{lem:normv} can be applied. 
Using Lemma~\ref{lem:normv}, a sufficient condition for this lemma
is to enforce
\begin{equation}
v(0) - t_{max}e^{\gamma_0}(e^{v_{max}} - 1) - t_{max}C_1\exp\left(C_2\right)
\geq v_{min}. 
\end{equation}
Rearranging terms and expanding $C_2$ we find
\begin{equation}
v(0) - t_{max}e^{\gamma_0}(e^{v_{max}} - 1) - v_{min} \geq
t_{max}C_1\exp\left(\beta\|T\|_1 t_{max}\right),
\end{equation}
where solving for $\|T\|_1$ completes the proof. 
\hfill $\blacksquare$

With Assumption~\ref{as:v} in place, 
Lemmas~\ref{lem:vmax} and~\ref{lem:vmin} jointly provide that 
$v(t) \in \big[v_{min}, v_{max}\big]$ for all~$t \in [0, t_{max}]$,
and the next two lemmas rely on this point.

\subsection{Enforcing $\dot{v} \leq d_{max}$}
In this section we derive a condition on $\|T\|_1$ that implies 
$\dot{v}(t) \leq d_{max}$ for all $t \in [0, t_{max}]$. 

\begin{lemma}
Suppose that $\|T\|_1$ satisfies the bound in Lemma~\ref{lem:vmax}. 
If $\|T\|_1$ also satisfies
\begin{equation}
\|T\|_1 \leq \frac{1}{\beta t_{max}}\ln\left(\frac{d_{max} + e^{\gamma_0}(e^{v_{min}} - 1)}{C_1}\right)
\end{equation}
then $\dot{v}(t) \leq d_{max}$ for all $t \in [0, t_{max}]$. 
\end{lemma}
\emph{Proof:} Similar to above, we use that
\begin{equation}
\dot{v} = -e^{\gamma_0}(e^v - 1) - \bone^Tw \leq -e^{\gamma_0}(e^{v_{min}} - 1)
+ \|w(t)\|_1
\end{equation}
to derive a sufficient condition. In particular, 
because~$T$ satisfies the bound in Lemma~\ref{lem:vmax} 
we may
apply Lemma~\ref{lem:normw}. Doing so, 
we wish to enforce
\begin{equation}
-e^{\gamma_0}(e^{v_{min}} - 1) + 
C_1\exp\left(C_2\right) \leq d_{max}. 
\end{equation}
Rearranging we find
\begin{equation}
C_1\exp\left(\beta \|T\|_1 t_{max}\right) \leq d_{max} + e^{\gamma_0}(e^{v_{min}} - 1), 
\end{equation}
where isolating $\|T\|_1$ completes the proof. \hfill $\blacksquare$

\subsection{Enforcing $\dot{v} \geq d_{min}$}
This section derives a sufficient condition on $\|T\|_1$ that implies
$\dot{v} \geq d_{min}$ for all $t \in [0, t_{max}]$. 

\begin{lemma}
Suppose that~$T$ obeys the bound in Lemma~\ref{lem:vmax}. 
If $\|T\|_1$ also obeys the bound
\begin{equation}
\|T\|_1 \leq \frac{1}{\beta t_{max}}\ln\left(\frac{|d_{min}| - e^{\gamma_0}(e^{v_{max}} - 1)}{C_1}\right),
\end{equation}
then $\dot{v}(t) \geq d_{min}$ for all $t \in [0, t_{max}]$. 
\end{lemma}
\emph{Proof:}
Using Equation~\eqref{eq:v}, a sufficient condition for $\dot{v}(t) \geq d_{min}$ is
$d_{min} \leq -e^{\gamma_0}(e^v - 1) - \bone^Tw$. 
Rearranging, we wish to enforce
$e^{\gamma_0}(e^v - 1) + \bone^Tw \leq |d_{min}|$, 
where the absolute value comes from the fact that $d_{min} < 0$
in Assumption~\ref{as:v}. 
A sufficient condition for doing so is 
$e^{\gamma_0}(e^v - 1) + \|w(t)\|_1 \leq |d_{min}|$.
By hypothesis,~$T$ satisfies the bound in Lemma~\ref{lem:vmax}
and, as as result,~$v(t) \leq v_{max}$ for all~$t \in [0, t_{max}]$. Then
we apply Lemma~\ref{lem:normw} to find the bound
\begin{equation}
e^{\gamma_0}(e^{v_{max}} - 1) + C_1\exp\left(\beta \|T\|_1t_{max}\right)
\leq |d_{min}|. 
\end{equation}
Solving for $\|T\|_1$ then gives the result. \hfill $\blacksquare$

\subsection{Overall Sustainability Bound}
Enforcing sustainability requires that~$v(t) \in [v_{min}, v_{max}]$
and that~$\dot{v}(t) \in [d_{min}, d_{max}]$ for all~$t \in [0, t_{max}]$.
Insuring satisfaction of these conditions can be done by satisfying
all four of the preceding bounds on~$\|T\|_1$. 
Combining the four preceding lemmas, we have the following theorem
that states a single unified sustainability criterion. 

\begin{theorem} \label{thm:main}
Suppose that Assumption~\ref{as:v} holds and define the constants
\begin{align}
\xi_1 &= e^{\gamma_0}(e^{v_{max}} - 1) \\
\xi_2 &= \frac{v(0) - t_{max}e^{\gamma_0}(e^{v_{max}} - 1) - v_{min}}{t_{max}} \\
\xi_3 &= d_{max} + e^{\gamma_0}(e^{v_{min}} - 1) \\
\xi_4 &= -d_{min} - e^{\gamma_0}(e^{v_{max}} - 1).
\end{align}
If, for all $i \in \{1, 2, 3, 4\}$, 
$\xi_i > C_1$,
then the system in Equation~\eqref{eq:vwdynamics} 
is sustainable
for any~$T$ satisfying
\begin{equation} \label{eq:sust_cond}
\|T\|_1 \leq \frac{1}{\beta t_{max}} \ln\left(\frac{\min\limits_{i \in \{1, 2, 3, 4\}} \xi_i}{C_1}\right).
\end{equation}
\hfill $\blacksquare$
\end{theorem}

The condition that $\xi_i > C_1$ is simply a feasibility condition; if it is
violated, then these bounds cannot be used to insure sustainability because
the logarithm will output a negative value. This feasibility condition can be used to generate conditions that must be satisfied by the parameters $v_{min},v_{max},d_{min},d_{max}$, and $t_{max}$;
we avoid a lengthy exposition on this subject because the many parameters
in this model generate a vast number of possible relationships that
one could explore. Instead, we comment on one relationship in particular.
Elementary operations show that the condition $\xi_{1}>C_{1}$ implies that $\frac{1}{\sum_{i=1}^{n} b_{i}\alpha_{i} }>t_{max}$. This relationship points to the importance of agent sensitivity, captured by $b_{i}$, when ensuring sustainability. We expand on the relationship between sustainability and sensitivity in the next section, and we further discuss the implications
of sustainability in that context. 

\section{Simulation and Discussion} \label{sec:sim}
In this section we explore the implications of the above results for the behavior of social-ecological systems. First, the behavior of the model is examined under different $\theta$ parameters in order to demonstrate the differences between societies that place more weight on social information and societies that place more weight on ecological information; these cases will be referred to as \emph{pro-social} and \emph{pro-ecological} societies, respectively. 
We will benchmark these cases against societies which weight social and
ecological information comparably, which we refer to as the
\emph{equal} case. 
Second, the sustainability bounds are discussed with respect to the parameter $b_{i}$, demonstrating how the sensitivity of populations relates to the sustainability
of their consumptions.  

\subsection{Pro-Social and Pro-Ecological Communities}
The stability proof presented in Theorem \ref{thm:stability} depends on Assumption \ref{as:weights}, i.e., that $\theta_i \geq \sum_{\substack{k = 1 \\ k \neq i}}^{n} \omega_{ki}\theta_k$ for all~$i$. 
If 
the collection
$\{\theta_{i}\}_{i \in \{1, \ldots, n\}}$ satisfies Assumption \ref{as:weights}, then $\{\delta\theta_{i}\}_{i \in \{1, \ldots, n\}}$ for $\delta\geq0$ also does. While the relationship between the individual $\theta_{i}$'s is fixed by Assumption \ref{as:weights}, the scaling factor $\delta$ allows a 
variety of systems to be described based on a single collection $\{\theta_{i}\}_{i \in \{1, \ldots, n\}}$ that satisfies Assumption~\ref{as:weights}. 
With such a collection,
and using $\theta_{i}=\frac{\nu_{i}}{\alpha_{i}}$ and $\nu_{i}+\alpha_{i}=1$, we set
\begin{equation}
\alpha_{i}=\frac{1}{1+\delta\theta_{i}} \ \ \text{ and } \ \ \nu_{i}=\frac{\delta\theta_{i}}{1+\delta\theta_{i}}
\end{equation}
to find new sets of admissible parameters. 
If $\delta$ is chosen so that $\delta\theta_{i}\gg 1$, 
then~$\alpha_i \approx 0$
and~$\nu_{i}\approx 1$, which we refer to as the \emph{pro-social} condition. 
If~$\delta\theta_{i}\ll 1$, then~$\nu_{i}\approx 0$ and~$\alpha_i \approx 1$, which we refer to as the \emph{pro-ecological} condition.
We refer to the case of~$\alpha_i \approx \nu_i$ as
the \emph{equal} condition. 

\begin{figure}
\centering
\includegraphics[width=.65\columnwidth]{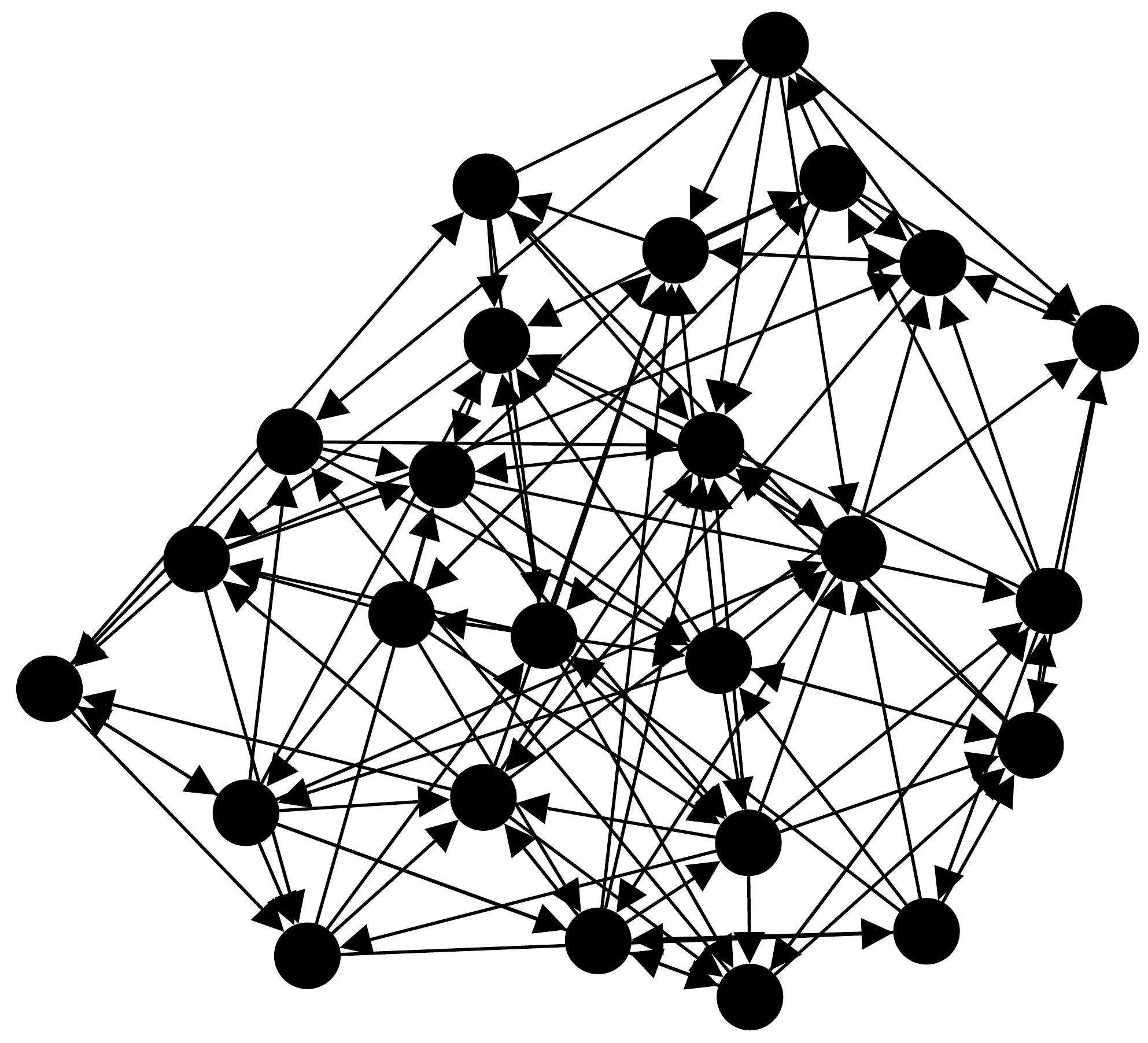}
\caption{The $25$-node, $114$-edge random graph used for simulation visualized with Gephi~\cite{ICWSM09154}. 
Each node in this graph represents an agent and each edge represents
that two agents interact. 
}
\label{fig:graph}
\end{figure}
\begin{figure}
\includegraphics[width=\columnwidth]{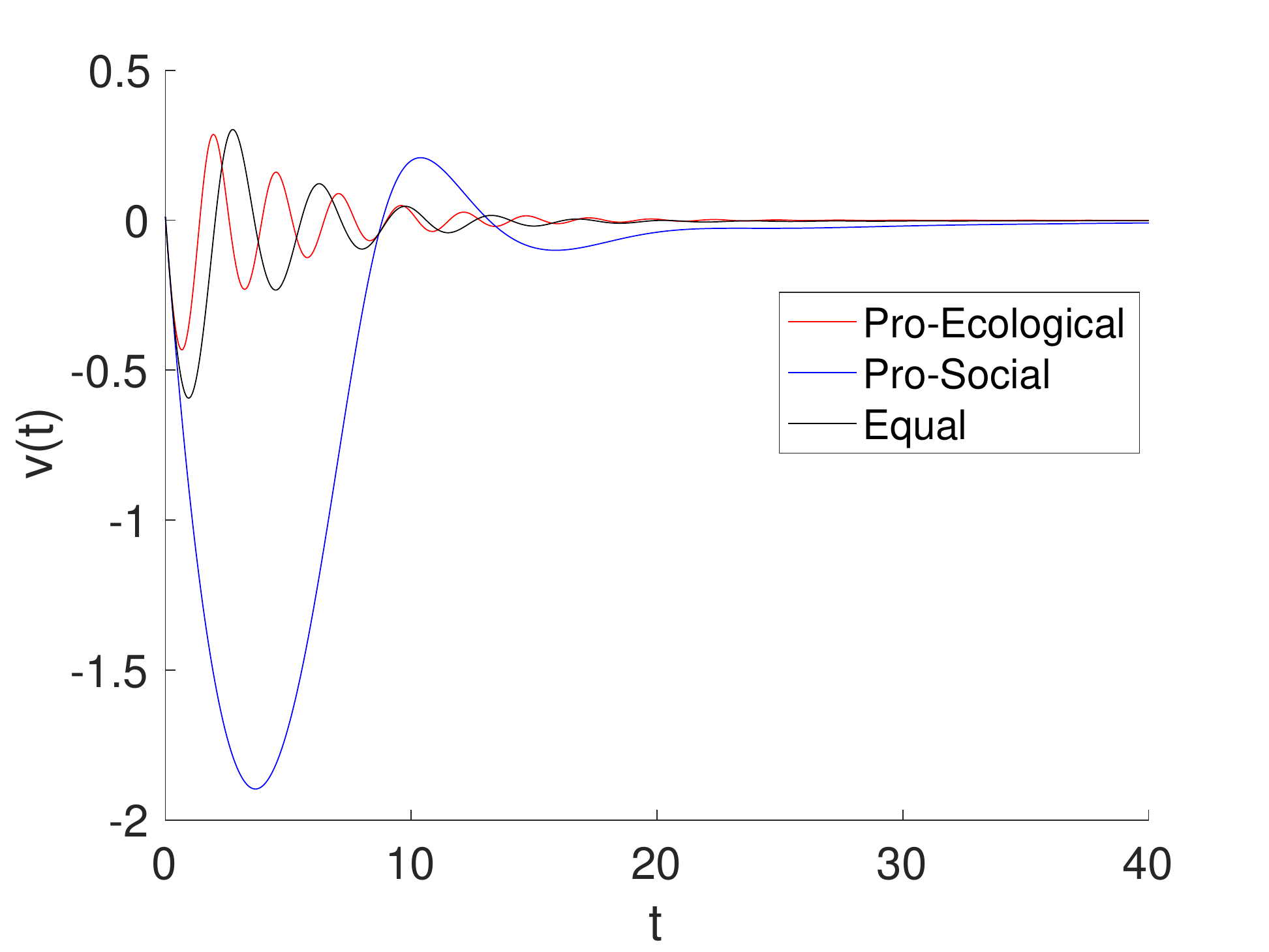}
\caption{
The behavior of the resource stock~$v(t)$ for the pro-ecological,
pro-social, and equal cases when simulated on the network
in Figure~\ref{fig:graph}. 
The pro-ecological and equal cases show similar oscillatory
behaviors while the pro-social cases shows lower-frequency oscillations
of substantially larger magnitude. 
} 
\label{fig:stab_v}
\end{figure}

To understand the impact of the relative sizes of~$\alpha_{i}$ 
and~$\nu_{i}$ on system behavior, the system in Equation~\eqref{eq:vwdynamics} was simulated on the $25$-node randomly
generated 
graph shown in Figure \ref{fig:graph}. Each edge was given uniform weight. The choice of~$\theta= [0.1826~~ 0.3296~~    0.2313  ~~  0.3454  ~~  0.1987  ~~  0.1923  ~~  0.1642  ~~ 0.1989  ~~  0.1182~~ 0.2198  ~~  0.1124 ~~  0.0734  ~~  0.1592 \\  0.3608  ~~  0.1913  ~~  0.1810   ~~ 0.2098  ~~  0.1206~~ 0.3210  ~~  0.0606  ~~  0.0597  ~~ 0.1302  ~~  0.0808  ~~  0.1336  ~~  0.1638]^{T}$ insures
satisfaction of Assumption~\ref{as:weights}. 
Three cases were run to elucidate the impacts
of varying~$\alpha_i$ and~$\nu_i$. In terms of their respective
averages,~$\bar{\alpha}$ and~$\bar{\nu}$, these cases are: 
\begin{enumerate}[i.]
\item The pro-ecological case with $\delta=0.1$,
which gives $\bar{\alpha}=0.9822$ and $\bar{\nu}=0.0196$
\item The pro-social case with $\delta=10$, which gives ${\bar{\alpha}=0.0644}$ and ${\bar{\nu}=0.9356}$ 
\item The \emph{equal} case with $\delta=\bar{\theta}=0.1816$, which gives $\bar{\alpha}=0.5261$ and $\bar{\nu}=0.4739$.
\end{enumerate}
The parameters $b_{i}$ and $\rho_{i}$ were chosen uniformly at random from $[0,1]$.

\begin{figure*}
\includegraphics[width=.33\columnwidth]{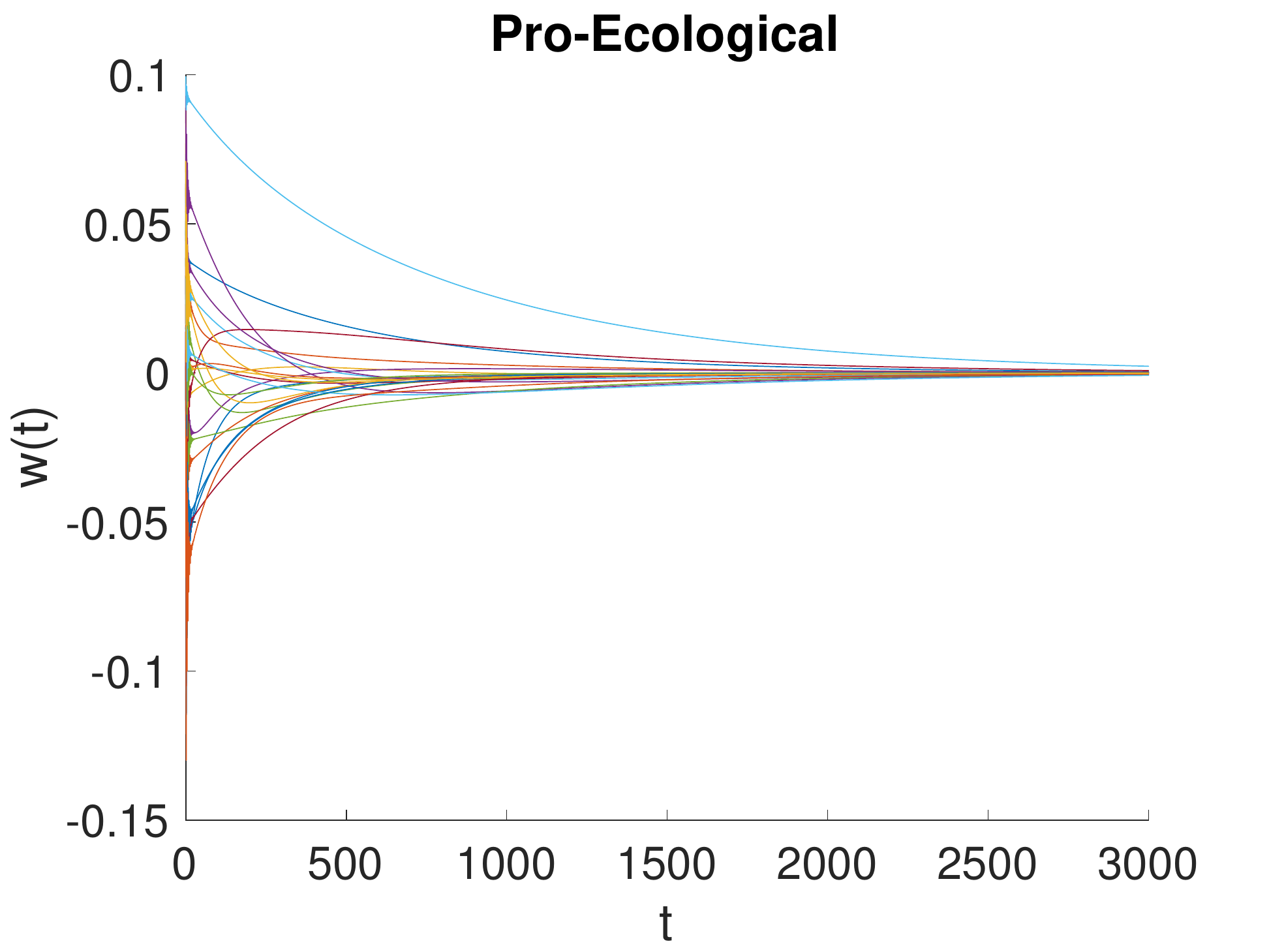}
\includegraphics[width=.33\columnwidth]{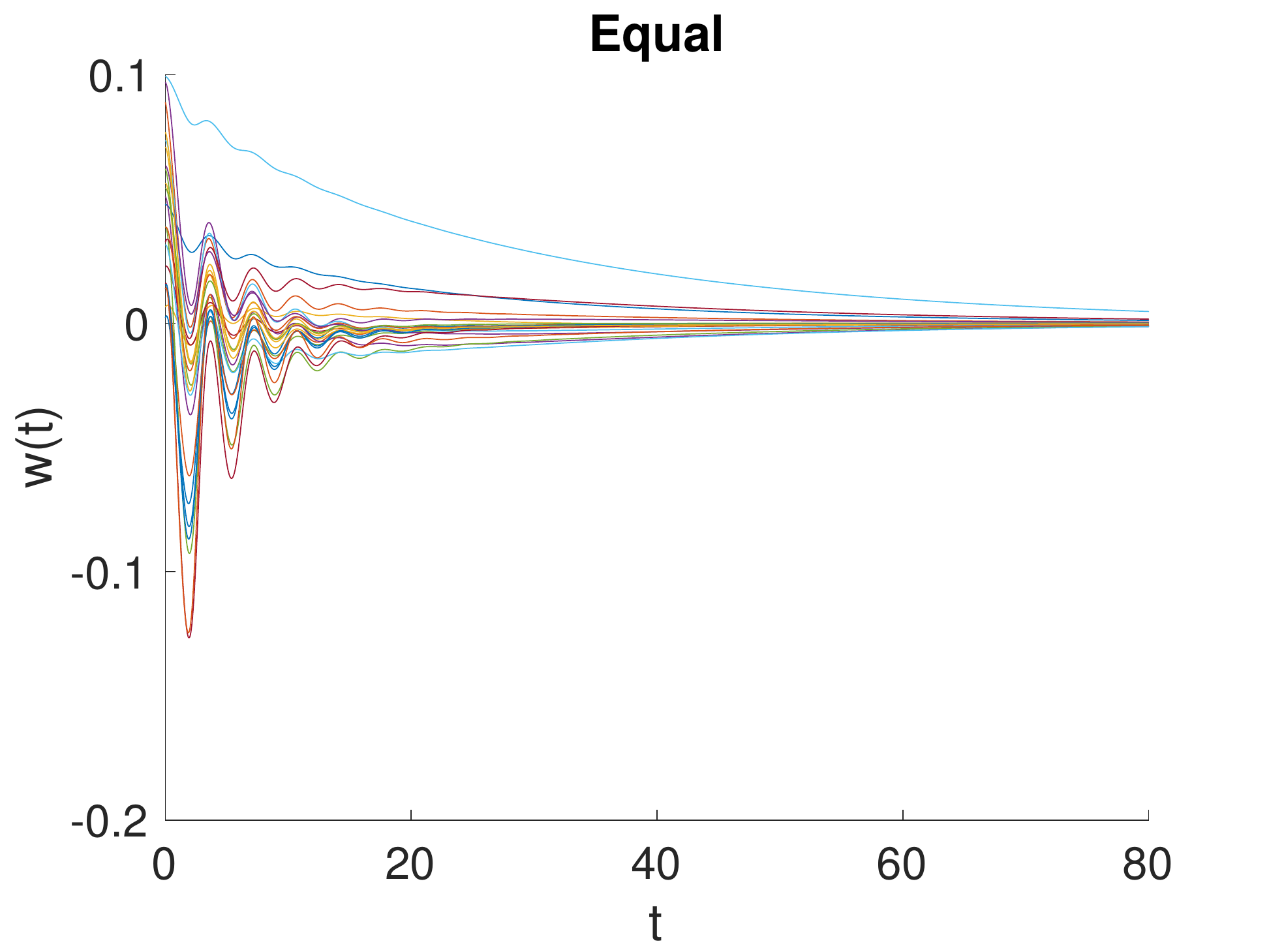}
\includegraphics[width=.33\columnwidth]{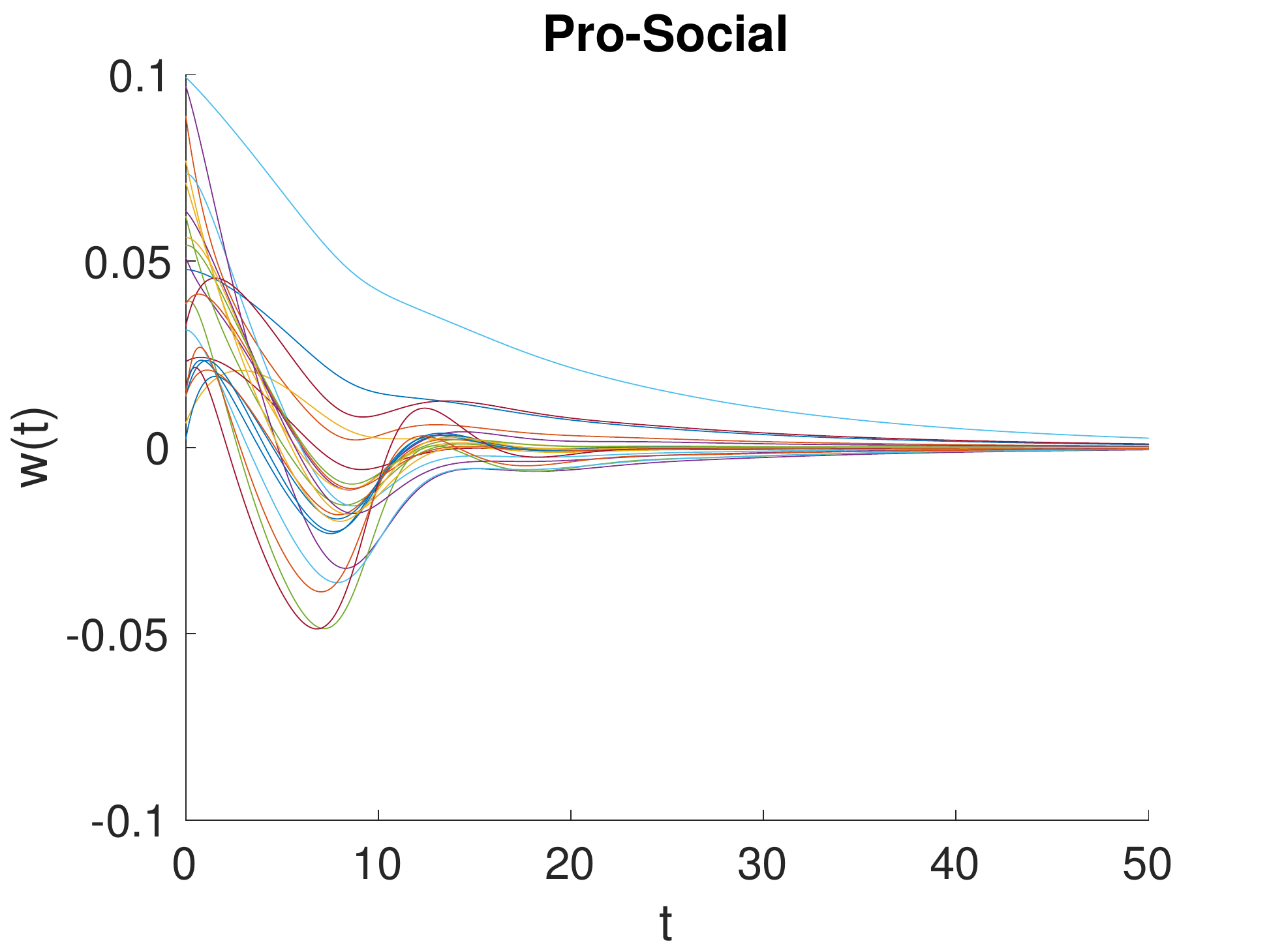}
\caption{The individual consumptions, $w_i(t)$, $i \in \{1, \ldots, 25\}$, for 
the aforementioned $25$-node graph under 
the pro-social, pro-ecological, and equal conditions. 
It takes approximately~$2000$ time steps for the pro-ecological network to converge, compared with~$100$ and~$50$ time steps for the equal and pro-social networks, respectively. }
\label{fig:stab_w}
\end{figure*}
The resource behavior for all three runs is shown in Figure \ref{fig:stab_v} and the individual consumption behavior  for all three runs is shown in Figure \ref{fig:stab_w}; 
each plot is representative of the behavior observed on a variety of graph topologies over the course of many simulations. 
We see that both the resource stock and consumption values 
tend to fluctuate to varying degrees, depicted by the oscillations in
these trajectories. Such oscillations may be undesirable in social-ecological systems for a variety of reasons. For instance, it has been observed that grasslands and forests exhibit an increased tendency to be eradicated by invading plant species in the presence of frequent disturbances and repeated exposure to extreme environmental conditions \cite{davis2000fluctuating, davis2001experimental}. 
In addition, boom-bust cycles in human-natural systems may not only lead to undesirable development patterns in the linked communities \cite{rodrigues2009boom}, but also unwanted price fluctuations by means of substantial supply shocks \cite{czech2013supply}.  

Among the three shown simulation runs, 
Figure~\ref{fig:stab_v} shows that the pro-social system exhibits relatively fewer
total  oscillations than the other systems, though the magnitude of oscillations
is larger.  
Note that in Figure~\ref{fig:stab_v}, while the resource dynamics for the pro-ecological and equal societies display significantly more fluctuations than the pro-social society, they converge to the equilibrium faster. The notion of speed of convergence, and of return to the equilibrium, have been 
previously associated with the concept of resilience, a fundamental 
attribute of sustainable 
systems~\cite{patten1997logical,walker2004resilience}. We then see from 
Figure~\ref{fig:stab_v} that while the pro-social resource dynamics have desirably fewer fluctuations, the pro-ecological dynamics have a desirably higher rate of convergence with smaller magnitudes of oscillations. 
The situation is reversed if we observe the trajectories of the individual consumptions in Figure~\ref{fig:stab_w}. Here, the pro-ecological society exhibits fewer oscillations but with slower convergence, while the pro-social society exhibits faster convergence but with more oscillations. 

There thus exist two trade-offs.
The first is between fluctuations and rate of convergence, with
faster convergence coming at the cost of more fluctuations in the resource,
and reduced fluctuations coming at the cost of slower convergence to
an equilibrium value. 
The second tradeoff is between desirable behaviors in the resource
and desirable behaviors in agents' consumption levels, with reduced
fluctuations in each coming at the cost of increased fluctuations in
the other. 
These tradeoffs suggest that the equal society, which gives a balanced preference to both ecological and social information, is more favorable than extreme societal configurations that weigh one source of information much more heavily 
than the other. From a policy perspective, 
promoting the equal society would entail measures that encourage society to give equal consideration to the state of the resource and the behavior of neighboring agents while determining individual consumptions. 
Examples of such measures include strategic information dispersion, targeted advertising, awareness campaigns, and manipulating visibility of key 
variables~\cite{steg2009encouraging}.


\subsection{System Sensitivity}
The conditions required by Theorem~\ref{thm:stability} did not restrict the parameters $\{b_{i}\}_{i \in \{1, \ldots, n\}}$, which factors out of the consumption dynamics 
for~$\dot{w}$ and thus essentially acts as a gain on agents' responsiveness to network stimuli. That is, the parameter~$b_i$ captures agent~$i$'s sensitivity to the social and ecological information they receive. However, the sustainability bounds derived in Section~\ref{sec:sust} depend upon $\{b_{i}\}_{i \in \{1, \ldots, n\}}$ through the parameter $\beta$ and the term $\sum_{i=1}^{n} b_{i} \alpha_{i}$ in the definition of $C_{1}$. Furthermore, as discussed previously, for $ \xi_{1}>C_{1}$ it must hold that $\frac{1}{\sum_{i=1}^{n} b_{i}\alpha_{i} }>t_{max}$, implying that agents' sensitivities
affect bounds on the time horizon over which sustainability can be guaranteed. 

We consider the implications of this bound in the same $25$-node network 
shown in Figure~\ref{fig:graph} under the equal parameter regime. We 
consider three uniform cases: $b_{i} = 0.025$, $b_i = 0.005$, or~$b_i = 
0.001$ for all~$i$, and 
changing~$b_i$ changes the time horizon over which a system 
can be shown to be sustainable. We consider a time horizon of $t_{max}=1$ 
when $b_{i}=0.025$, a time horizon of $t_{max}=2$ when $b_{i}
=0.005$, and a time horizon of $t_{max}=3$ when $b_{i}=0.001$. 
Note that these times are normalized with respect to the intrinsic growth 
rate of the resource, and thus sustainability for~$t_{max} = 3$ pertains to 
the system behavior across an interval that is~$3$ times the intrinsic 
time constant of the resource~\cite{perman2003natural}.

\begin{figure}
\makebox[\columnwidth][c]{\includegraphics[width=1\columnwidth]{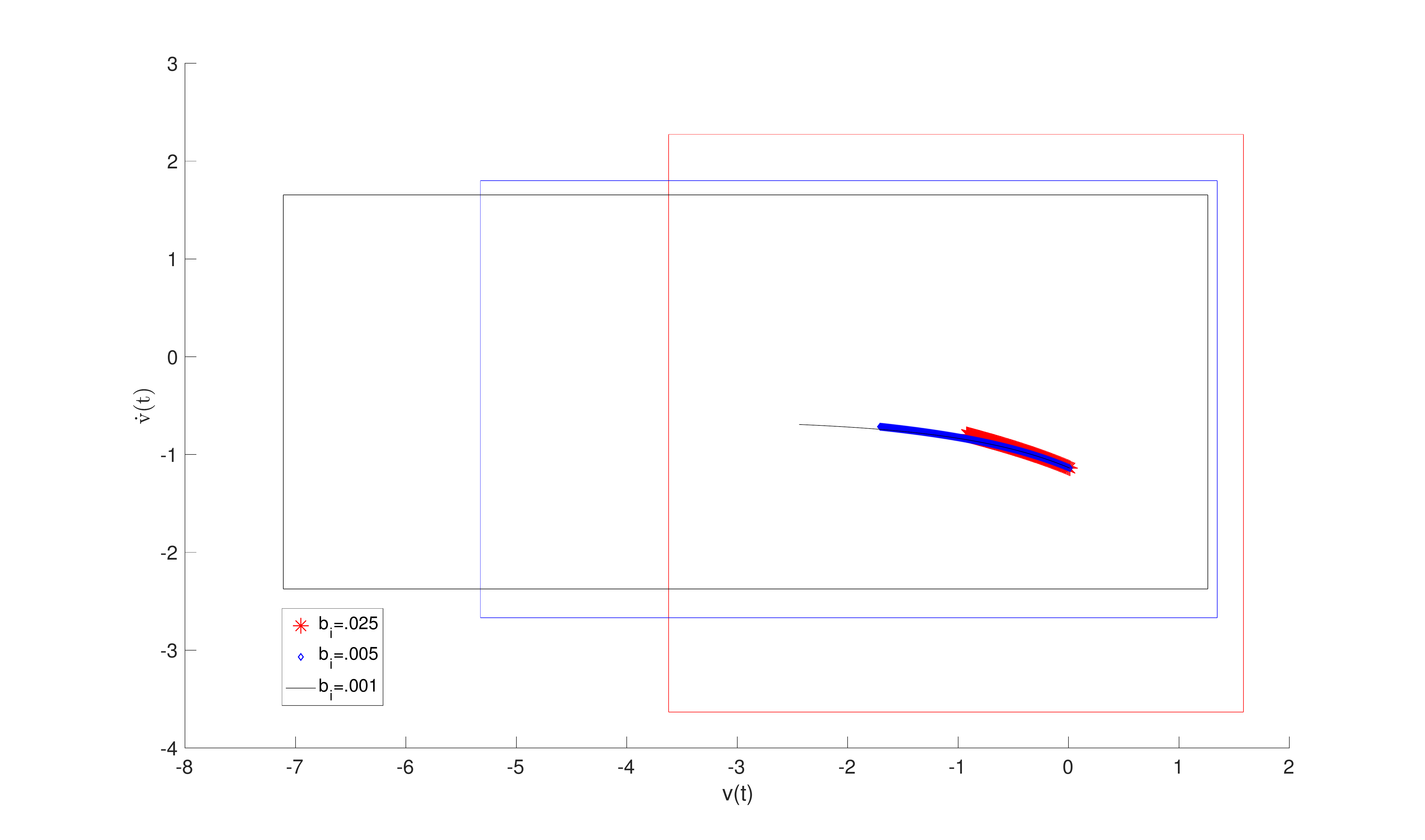}}
\caption{The resource, $v(t)$, versus  $\dot{v}(t)$ from $t=0$ to $t=t_{max}$ 
with~$b_{i}=0.025$ for all~$i$, $b_{i}=0.005$~for all~$i$, 
and~$b_{i}=0.001$~for all $i$ on the network in Figure~\ref{fig:graph} 
with~$\alpha_{i}=\nu_{i}$~for all $i$. 
The red, blue, and black boxes are the minimal sustainability window for $b_{i}=0.025$,
$b_{i}=0.005$, and $b_{i}=0.001$, respectively. As $b_{i}$ decreases, the trajectory 
elongates, reflecting both 
a larger $t_{max}$ and the fact that 
 the minimal sustainability window more tightly bounds 
the trajectories. 
} \label{fig:phase_plot}
\end{figure}

The bounds in Equation~\eqref{eq:sust_cond} determine sustainability with respect to choices 
of~$v_{min},v_{max},d_{min}$, and~$d_{max}$, which would be determined by the ecological context that is being modeled. For this simulation, 
we consider the minimal $v_{max}$ and $d_{max}$ as well as the maximal $v_{min}$ and $d_{min}$ for which the system is sustainable, which we refer to as the 
\emph{minimal sustainability window}. 
These values can be derived from elementary operations on the condition that $\xi_{i}=e^{\|T\|_1 \beta t_{max}}C_{1}$ for $i\in\{1,2,3,4\}$, which follows from attaining equality
in the bound in Theorem~\ref{thm:main}.  
We 
plot the evolution of the resource from $t=0$ to $t=t_{max}$ 
in Figure~\ref{fig:phase_plot}, where we see that that it indeed
remains within the minimal sustainability window. 
In Figure~\ref{fig:phase_plot}, as $b_{i}$ decreases the minimal sustainability window more tightly bounds the resource trajectory, 
suggesting that the resulting system is sustainable with respect to a larger set of choices of~$v_{min},v_{max},d_{min}$, and~$d_{max}$.  
This also suggests that a decrease in agent sensitivity results in longer 
time horizons for which the system can be made sustainable.

However, sensitive systems have been observed to be more beneficial in certain settings and with respect to a particular set of indicators. For instance, Consumer Affect (intensity of reaction to stimuli~\cite{burke2006contemporary}) has been found to have a profound effect on achieving sustainable consumer behavior~\cite{steg2009encouraging}. The ability of agents to adapt quickly to resource fluctuations has also been found to be an important aspect of sustainable fishing communities~\cite{allison2001livelihoods}. We find in our model that low sensitivity can indeed be associated with higher fluctuations in the resource. 

\begin{figure}
\includegraphics[width=\columnwidth]{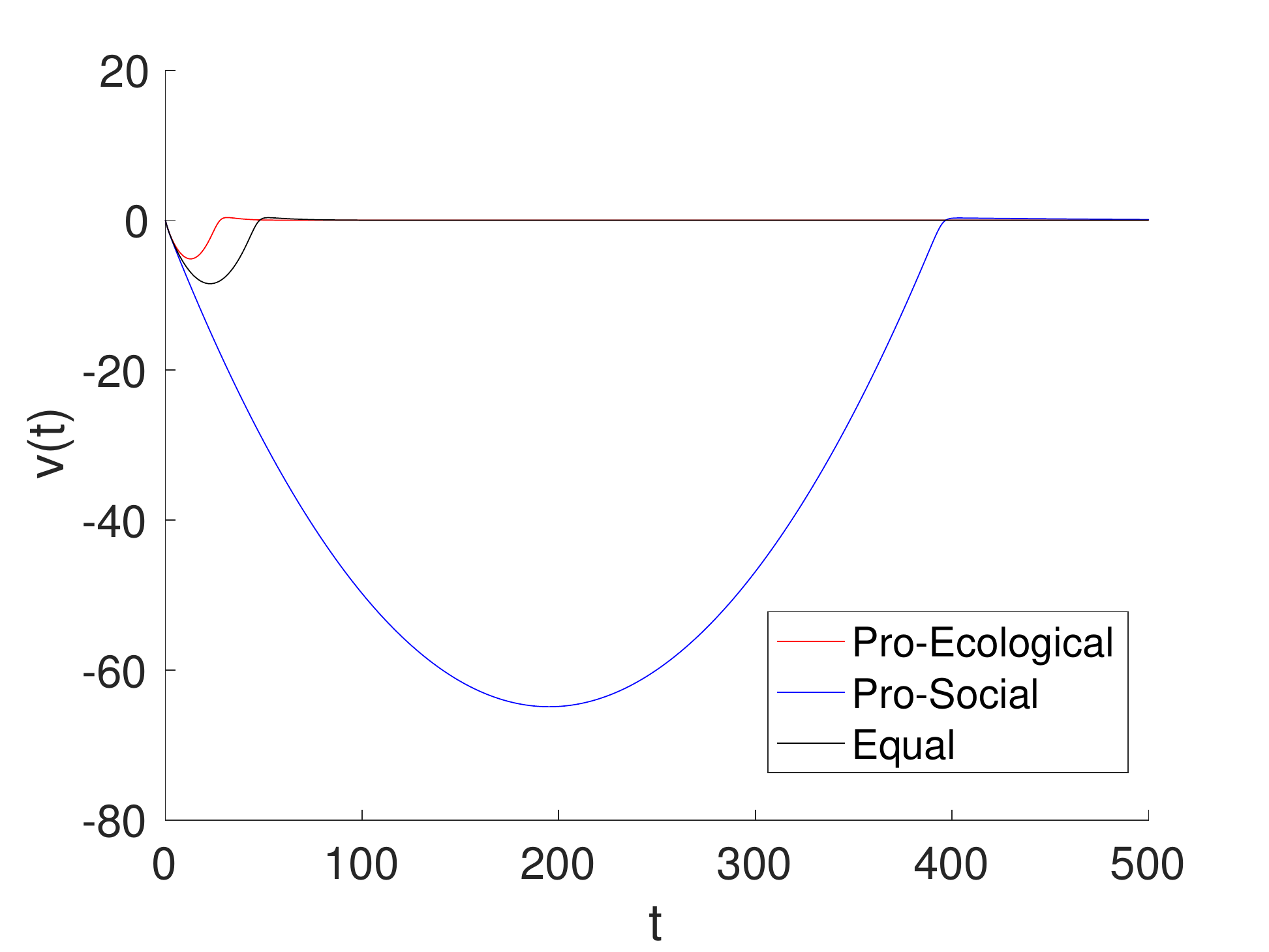}
\caption{
The behavior of the resource stock, $v(t)$, 
for the $25$-node graph in Figure~\ref{fig:graph} under three parameter conditions on $\alpha$ and $
\nu$ when $b_{i}=0.005$~for all $i$. 
Here the pro-social condition has oscillations that are much larger than the other two cases. 
}
\label{fig:low_b}
\end{figure}

However, as highlighted in the preceding sections, social-ecological systems are complex in nature, where a single factor alone cannot be associated with a particular outcome. Instead, system parameters and variables often act in combination to produce a certain outcome. 
To further illustrate this point, 
Figure~\ref{fig:low_b} shows the resource trajectory for a relatively less-sensitive society. 
The simulation was run for the same $25$-node graph as above,
with $b_{i}=0.005$ and for all configurations of information preference used above. The pro-social society has a very large fluctuation in the quantity of the resource, reaching\footnote{The negative value here is due to the 
logarithmic coordinate transformation to arrive at the
system in~$(v,w)$-coordinates.}  
$v\approx-30$,
which is roughly~$6$ times the size of fluctuations in the equal society and nearly~$10$ times that of the pro-ecological one. This suggests that insufficient weighting of ecological information and low sensitivity, together, lead to large, potentially harmful fluctuations in the resource. This point has been
observed previously in~\cite{uribe2012community}, where more sensitivity
to ecological information 
has been observed to contribute to adopting more sustainable 
behavior patterns.

The above suggests another trade-off enacted by varying the sensitivity of the agents. On one hand, higher sensitivity leads to loose bounds and small time horizons for the sustainability criterion. On the other hand, low sensitivity leads to larger oscillations in the resource stock (especially in combination with low preference for ecological information). The model hence suggests, from a policy perspective, controlled measures to vary the responsiveness of the consuming population. These measures may either be of an informational or structural nature (see~\cite{steg2009encouraging} and included references for possible examples of such measures), with the end goal always being
improving sustainability outcomes. 




\section{Conclusion}\label{sec:conc}
In this paper we have presented a mathematical criterion for sustainability grounded in the literature on 
sustainable development. 
The presented definition of sustainability extends the notion of stability which is often applied to ecological 
systems. While stability pertains to the asymptotic behavior of the system, sustainability concerns itself with 
behavior over a finite time horizon where transient behavior often dominates. 
We find that, in the chosen model of natural resource consumption, the requirements for sustainability are not 
captured by stability alone. In particular, we observe that agent sensitivity, while having no contribution to 
determining stability, has a major role in sustainability of the same system. Our simulations of sustainable networks 
for the model uncover tradeoffs in system behavior that do not appear when investigating 
system stability alone. 
These findings may be translated to guidelines for effective policy-making in natural resource systems. This 
study thus serves as one step towards a dynamical systems 
theory for sustainability of social-ecological systems.     
 
\appendix
Toward proving Theorem~\ref{thm:stability}, 
we begin with the following lemma.

\begin{lemma} \label{lem:semi}
For all $x \in \R^n$ 
\begin{equation} \label{eq:quadformTVA} 
x^T(T^TVA^{-1})x \geq 0. 
\end{equation}
\end{lemma}
\emph{Proof:} Define $\theta_i = \frac{\nu_i}{\alpha_i}$ and the matrix
$\Theta = \textnormal{diag}(\theta_1, \ldots, \theta_n)$. The matrix $T^TVA^{-1}$ is not
symmetric in general, but only its symmetric part will
contribute to the quadratic form above. 
Noting that $T^TVA^{-1} = T^T\Theta$, we find
\begin{equation} \label{eq:l1runback}
x^T(T^T\Theta)x = \frac{1}{2}x^T(T^T\Theta + \Theta T)x,
\end{equation}
where we have used the fact that $\Theta = \Theta^T$; $T \neq T^T$ in general
because we have not assumed symmetric weights. Expanding, we find
\begin{equation}
T^T\Theta + \Theta T  = 
\left(\begin{array}{ccc} 
2\theta_1 & \cdots & -\omega_{n1}\theta_n - \omega_{1n}\theta_1 \\
-\omega_{12}\theta_1 - \omega_{21}\theta_2 & 
-\omega_{n2}\theta_n - \omega_{2n} \theta_2 \\
\vdots & \ddots & \vdots \\
-\omega_{1n}\theta_1 - \omega_{n1}\theta_n & \cdots & 2\theta_n
\end{array}\right).
\end{equation}

This matrix has real eigenvalues
because it is symmetric. Applying Gershgorin's circle theorem to column~$i$ of 
$T^T\Theta + \Theta T$, we see that non-negativity of these eigenvalues
requires
\begin{equation}
2\theta_i \geq \sum_{\substack{j=1 \\ j \neq i}}^{n} \omega_{ij}\theta_i
+ \sum_{\substack{k=1 \\ k \neq i}}^{n} \omega_{ki}\theta_k 
= \theta_i + \sum_{\substack{k=1 \\ k \neq i}}^{n} \omega_{ki}\theta_k,
\end{equation}
which follows from the fact that 
$\sum_{\substack{j= 1 \\ j \neq i}}^{n} \omega_{ij} = 1$
under the model in this problem. 

This condition is then equivalent to requiring
\begin{equation}
\theta_i \geq \sum_{\substack{k=1 \\ k \neq i}}^{n} \omega_{ki}\theta_k,
\end{equation}
which holds under Assumption~\ref{as:weights}. Then the quadratic
form in Equation~\eqref{eq:quadformTVA} is non-negative. 
\hfill $\blacksquare$

Next, we show that $TVA^{-1}$ has a non-trivial nullspace.

\begin{lemma} \label{lem:onenull}
Let $\bone$ denote the vector of all $1$'s in $\R^n$. Then 
$\bone$ is in the nullspace of $T^T\Theta + \Theta T$. 
\end{lemma}
\emph{Proof:} Using
$\sum_{\substack{j= 1 \\ j \neq i}}^{n} \omega_{ij} = 1$,
we find
%
\begin{align}
\bone^T\big(T^T\Theta + \Theta T\big)\bone &=
\sum_{i=1}^{n} \theta_i - \sum_{i=1}^{n} \sum_{\substack{j = 1 \\ j \neq i}}^{n} \omega_{ji} \theta_j \\
&= \sum_{i=1}^{n} \theta_i - \sum_{j=1}^{n} \theta_{j} \sum_{\substack{i = 1 \\ i \neq j}}^{n} \omega_{ji} = 0, 
\end{align}
where we have used that
$\sum_{\substack{i=1 \\ i \neq j}}^{n} \omega_{ji} = 1$.
We conclude by using the fact that, for a symmetric and positive
semi-definite matrix~$M$, $x^TMx = 0$ if and only if
$x \in \mathcal{N}(M)$ (see Fact 8.15.2 in~\cite{bernstein09}). 
Then $\bone \in \mathcal{N}\big(T^T\Theta + \Theta T\big)$. 
\hfill $\blacksquare$

Lemma~\ref{lem:onenull} 
shows that $T^T\Theta + \Theta T$ is singular. In particular, it has
a nullspace of dimension at least $1$. In the next lemma we show that
the dimension of its nullspace is exactly $1$. 

\begin{lemma} \label{lem:rank}
The matrix $T^T\Theta + \Theta T$  has rank $n-1$.
\end{lemma}
\emph{Proof:} Strong connectivity, from Assumption~\ref{as:conn}, 
makes $T$ irreducible, so that $T^T$ is as well. Multiplying
$T^T \Theta$ only gives positive weights to the entries of $T^T$, and
$T^T \Theta$ is therefore also irreducible. Because all diagonal entries
are positive and all off-diagonal entries are non-positive, 
the sum $T^T \Theta + \Theta T$ is then irreducible as well. 
$T^T\Theta + \Theta T$ is also symmetric, and, by Lemma~\ref{lem:semi},
it is positive semidefinite. By Exercise~4.15 in~\cite{berman94}, it is
an~$M$-matrix. Theorem~4.16 in~\cite{berman94} shows that singular,
irreducible~$M$-matrices have rank~$n-1$. 
\hfill $\blacksquare$

Because $T^T\Theta + \Theta T$ has a one-dimensional nullspace containing
$\bone$, we see that
\begin{equation}
\mathcal{N}\big(T^T\Theta + \Theta T\big) = \textnormal{span}\{\bone\}.
\end{equation}

We now proceed with the main stability proof.

\emph{Proof of Theorem~\ref{thm:stability}:}
We use the Lyapunov function 
\begin{equation}
V(v, w) = e^v - v - 1 + \frac{1}{2}e^{-\gamma_0}w^T(AB)^{-1}w,
\end{equation}
where $(AB)^{-1}$ is well-defined because $A$ and $B$ are diagonal matrices
with positive diagonal elements. Differentiating $V$, we find that
\begin{align}
\dot{V} &= e^v\dot{v} - \dot{v} + \dot{w}^Te^{-\gamma_0}(AB)^{-1}w \\
        &= (e^v - 1)(-e^{\gamma_0}(e^v - 1) - \bone^Tw) + e^{\gamma_0}(e^v - 1)\bone^TABe^{-\gamma_0}(AB)^{-1}w - w^TT^TVBe^{-\gamma_0}(AB)^{-1}w,
\end{align}
where we have used that $A$ and $B$ are symmetric. 
Gathering like terms, we find that
\begin{equation}
\dot{V} = -e^{-\gamma_0}(e^v - 1)^2 - e^{-\gamma_0}w^TT^TVA^{-1}w. 
\end{equation}
The first term is manifestly negative for all $v \neq 0$. 
For $v = 0$, one has 
\begin{align} \label{eq:wconlater}
\dot{V} &= - e^{-\gamma_0}w^TT^TVA^{-1}w \\ &= -e^{-\gamma_0}w^T(T^T\Theta + \Theta T)w
\leq 0. 
\end{align}

Because $T^T\Theta + \Theta T$ is not full rank, we must apply LaSalle's invariance
principle to show asymptotic convergence. To do so, we wish to show that the set
\begin{equation}
S := \{(v, w) \mid \dot{V}(v, w) \equiv 0\}
\end{equation}
contains only the trivial trajectory $(v, w) \equiv 0$. We first note that
$\dot{V} = 0$ requires $v = 0$ for all time. 
Setting~$v = 0$ and~$\dot{v} = 0$,
we use the $v$ dynamics of this system to find
\begin{equation}
\dot{v} = 0 = -\bone^Tw,
\end{equation}
and thus we require $\bone^Tw = 0$ for~$w$ to remain in $S$, which we
write as $w(t) \in \textnormal{span}\{\bone\}^{\perp}$ for all~$t$.  
Next, from Lemma~\ref{lem:onenull}, we see that
\begin{equation}
w^T\big(T^T\Theta + \Theta T\big)w = 0
\end{equation}
if and only if $w$ is in the nullspace of $T^T\Theta + \Theta T$. 
As a result, with $v = 0$, Equation~\eqref{eq:wconlater} shows
we remain in $S$ if and only if $w(t) \in \textnormal{span}\{\bone\}$ for all~$t$. 
Combined with the above, remaining in $S$ requires
\begin{equation}
w(t) \in \textnormal{span}\{\bone\} \cap \textnormal{span}\{\bone\}^{\perp}
\end{equation}
for all time, 
which necessarily implies that $w \equiv 0$. 
\hfill $\blacksquare$

\bibliographystyle{IEEEtran}{}
\bibliography{sources}

\begin{thebibliography}{10}
\providecommand{\url}[1]{#1}
\csname url@samestyle\endcsname
\providecommand{\newblock}{\relax}
\providecommand{\bibinfo}[2]{#2}
\providecommand{\BIBentrySTDinterwordspacing}{\spaceskip=0pt\relax}
\providecommand{\BIBentryALTinterwordstretchfactor}{4}
\providecommand{\BIBentryALTinterwordspacing}{\spaceskip=\fontdimen2\font plus
\BIBentryALTinterwordstretchfactor\fontdimen3\font minus
  \fontdimen4\font\relax}
\providecommand{\BIBforeignlanguage}[2]{{%
\expandafter\ifx\csname l@#1\endcsname\relax
\typeout{** WARNING: IEEEtran.bst: No hyphenation pattern has been}%
\typeout{** loaded for the language `#1'. Using the pattern for}%
\typeout{** the default language instead.}%
\else
\language=\csname l@#1\endcsname
\fi
#2}}
\providecommand{\BIBdecl}{\relax}
\BIBdecl

\bibitem{ostrom2015governing}
E.~Ostrom, \emph{Governing the commons}.\hskip 1em plus 0.5em minus 0.4em\relax
  Cambridge university press, 2015.

\bibitem{brundtland1987our}
G.~Brundtland \emph{et~al.}, ``Our common future: Report of the 1987 world
  commission on environment and development,'' \emph{United Nations, Oslo},
  vol.~1, p.~59, 1987.

\bibitem{cabezas2002towards}
H.~Cabezas and B.~D. Fath, ``Towards a theory of sustainable systems,''
  \emph{Fluid phase equilibria}, vol. 194, pp. 3--14, 2002.

\bibitem{chichilnisky1995green}
G.~Chichilnisky, G.~Heal, and A.~Beltratti, ``The green golden rule,''
  \emph{Economics Letters}, vol.~49, no.~2, pp. 175--179, 1995.

\bibitem{martinet2009defining}
V.~Martinet \emph{et~al.}, ``Defining sustainability objectives,'' in
  \emph{Atelier d'economie des ressources naturelles et de l'environnement de
  Montr{\'e}al}, 2009, pp. 29--p.

\bibitem{ludwig1997sustainability}
D.~Ludwig, B.~Walker, and C.~S. Holling, ``Sustainability, stability, and
  resilience,'' \emph{Conservation ecology}, vol.~1, no.~1, 1997.

\bibitem{patten1997logical}
B.~C. Patten and R.~Costanza, ``Logical interrelations between four
  sustainability parameters: stability, continuation, longevity, and health,''
  \emph{Ecosystem Health}, vol.~3, no.~3, pp. 136--142, 1997.

\bibitem{kinzig2014consumption}
A.~Kinzig and C.~Perrings, ``Consumption, stability, and sustainability in
  social-ecological systems,'' in \emph{Sustainable Consumption:
  Multi-disciplinary Perspectives In Honour of Professor Sir Partha
  Dasgupta}.\hskip 1em plus 0.5em minus 0.4em\relax Oxford University Press,
  2014, pp. 221--245.

\bibitem{wang2010emergence}
F.-Y. Wang, ``The emergence of intelligent enterprises: From cps to cpss,''
  \emph{IEEE Intelligent Systems}, vol.~25, no.~4, pp. 85--88, 2010.

\bibitem{manzoor2016game}
T.~Manzoor, E.~Rovenskaya, and A.~Muhammad, ``Game-theoretic insights into the
  role of environmentalism and social-ecological relevance: A cognitive model
  of resource consumption,'' \emph{Ecological modelling}, vol. 340, pp. 74--85,
  2016.

\bibitem{prell2011social}
C.~Prell and {\"O}.~Bodin, \emph{Social Networks and Natural Resource
  Management: Uncovering the social fabric of environmental governance}.\hskip
  1em plus 0.5em minus 0.4em\relax Cambridge University Press, 2011.

\bibitem{videras2013social}
J.~Videras, ``Social networks and the environment,'' \emph{Annu. Rev. Resour.
  Econ.}, vol.~5, no.~1, pp. 211--226, 2013.

\bibitem{videras2012influence}
J.~Videras, A.~L. Owen, E.~Conover, and S.~Wu, ``The influence of social
  relationships on pro-environment behaviors,'' \emph{Journal of Environmental
  Economics and Management}, vol.~63, no.~1, pp. 35--50, 2012.

\bibitem{prell2009stakeholder}
C.~Prell, K.~Hubacek, and M.~Reed, ``Stakeholder analysis and social network
  analysis in natural resource management,'' \emph{Society and Natural
  Resources}, vol.~22, no.~6, pp. 501--518, 2009.

\bibitem{crowe2007search}
J.~A. Crowe, ``In search of a happy medium: How the structure of
  interorganizational networks influence community economic development
  strategies,'' \emph{Social Networks}, vol.~29, no.~4, pp. 469--488, 2007.

\bibitem{ramirez2009impact}
S.~Ramirez-Sanchez and E.~Pinkerton, ``The impact of resource scarcity on
  bonding and bridging social capital: the case of fishers’
  information-sharing networks in loreto, bcs, mexico,'' \emph{Ecology and
  Society}, vol.~14, no.~1, 2009.

\bibitem{manzoor2018learning}
T.~Manzoor, E.~Rovenskaya, A.~Davydov, and A.~Muhammad, ``Learning through
  fictitious play in a game-theoretic model of natural resource consumption,''
  \emph{IEEE Control Systems Letters}, vol.~2, no.~1, pp. 163--168, 2018.

\bibitem{rockenbauch2017social}
T.~Rockenbauch and P.~Sakdapolrak, ``Social networks and the resilience of
  rural communities in the global south: a critical review and conceptual
  reflections,'' \emph{Ecology and Society}, vol.~22, no.~1, 2017.

\bibitem{manzoor2017structural}
T.~Manzoor, E.~Rovenskaya, and A.~Muhammad, ``Structural effects and
  aggregation in a social-network model of natural resource consumption,''
  \emph{IFAC-PapersOnLine}, vol.~50, no.~1, pp. 7675--7680, 2017.

\bibitem{ruf2018stability}
S.~F. Ruf, M.~T. Hale, T.~Manzoor, and A.~Muhammad, ``Stability of leaderless
  resource consumption networks (to appear),'' in \emph{Decision and Control,
  2018 57th IEEE Conference on}.\hskip 1em plus 0.5em minus 0.4em\relax IEEE,
  2018.

\bibitem{perman2003natural}
R.~Perman, Y.~Ma, J.~McGilvray, and M.~Common, \emph{Natural resource and
  environmental economics}.\hskip 1em plus 0.5em minus 0.4em\relax Pearson
  Education, 2003.

\bibitem{festinger1954theory}
L.~Festinger, ``A theory of social comparison processes,'' \emph{Human
  relations}, vol.~7, no.~2, pp. 117--140, 1954.

\bibitem{mosler2003integrating}
H.-J. Mosler and W.~M. Brucks, ``Integrating commons dilemma findings in a
  general dynamic model of cooperative behavior in resource crises,''
  \emph{European Journal of Social Psychology}, vol.~33, no.~1, pp. 119--133,
  2003.

\bibitem{pezzey2017economics}
J.~C. Pezzey and M.~A. Toman, \emph{The Economics of Sustainability}.\hskip 1em
  plus 0.5em minus 0.4em\relax Routledge, 2017.

\bibitem{valente2005sustainable}
S.~Valente, ``Sustainable development, renewable resources and technological
  progress,'' \emph{Environmental and Resource Economics}, vol.~30, no.~1, pp.
  115--125, 2005.

\bibitem{meadows2012limits}
D.~Meadows and J.~Randers, \emph{The limits to growth: the 30-year
  update}.\hskip 1em plus 0.5em minus 0.4em\relax Routledge, 2012.

\bibitem{ben2012cybernetics}
M.~Ben-Eli, ``The cybernetics of sustainability: definition and underlying
  principles,'' \emph{Enough for All forever: A Handbook for Learning about
  Sustainability, Murray J, Cawthorne G, Dey C and Andrew C (eds.). Champaign,
  IL, Common Ground Publishing: University of Illinois}, vol.~14, 2012.

\bibitem{bhatia2002stability}
N.~P. Bhatia and G.~P. Szeg{\"o}, \emph{Stability theory of dynamical
  systems}.\hskip 1em plus 0.5em minus 0.4em\relax Springer Science \& Business
  Media, 2002.

\bibitem{leine2010historical}
R.~I. Leine, ``The historical development of classical stability concepts:
  Lagrange, poisson and lyapunov stability,'' \emph{Nonlinear Dynamics},
  vol.~59, no. 1-2, p. 173, 2010.

\bibitem{holling1973resilience}
C.~S. Holling, ``Resilience and stability of ecological systems,'' \emph{Annual
  review of ecology and systematics}, vol.~4, no.~1, pp. 1--23, 1973.

\bibitem{loucks1997quantifying}
D.~P. Loucks, ``Quantifying trends in system sustainability,''
  \emph{Hydrological Sciences Journal}, vol.~42, no.~4, pp. 513--530, 1997.

\bibitem{ulph2014sustainable}
A.~Ulph and D.~Southerton, \emph{Sustainable consumption: Multi-disciplinary
  perspectives in honour of Professor Sir Partha Dasgupta}.\hskip 1em plus
  0.5em minus 0.4em\relax Oxford University Press, USA, 2014.

\bibitem{kharrazi2013quantifying}
A.~Kharrazi, E.~Rovenskaya, B.~D. Fath, M.~Yarime, and S.~Kraines,
  ``Quantifying the sustainability of economic resource networks: An ecological
  information-based approach,'' \emph{Ecological Economics}, vol.~90, pp.
  177--186, 2013.

\bibitem{costanza1995defining}
R.~Costanza and B.~C. Patten, ``Defining and predicting sustainability,''
  \emph{Ecological economics}, vol.~15, no.~3, pp. 193--196, 1995.

\bibitem{ames1997inequalities}
W.~F. Ames and B.~Pachpatte, \emph{Inequalities for differential and integral
  equations}.\hskip 1em plus 0.5em minus 0.4em\relax Elsevier, 1997, vol. 197.

\bibitem{ICWSM09154}
M.~Bastian, S.~Heymann, and M.~Jacomy, ``Gephi: An open source software for
  exploring and manipulating networks,'' 2009.

\bibitem{davis2000fluctuating}
M.~A. Davis, J.~P. Grime, and K.~Thompson, ``Fluctuating resources in plant
  communities: a general theory of invasibility,'' \emph{Journal of Ecology},
  vol.~88, no.~3, pp. 528--534, 2000.

\bibitem{davis2001experimental}
M.~A. Davis and M.~Pelsor, ``Experimental support for a resource-based
  mechanistic model of invasibility,'' \emph{Ecology letters}, vol.~4, no.~5,
  pp. 421--428, 2001.

\bibitem{rodrigues2009boom}
A.~S. Rodrigues, R.~M. Ewers, L.~Parry, C.~Souza, A.~Ver{\'\i}ssimo, and
  A.~Balmford, ``Boom-and-bust development patterns across the amazon
  deforestation frontier,'' \emph{Science}, vol. 324, no. 5933, pp. 1435--1437,
  2009.

\bibitem{czech2013supply}
B.~Czech, \emph{Supply shock: economic growth at the crossroads and the steady
  state solution}.\hskip 1em plus 0.5em minus 0.4em\relax New Society
  Publishers, 2013.

\bibitem{walker2004resilience}
B.~Walker, C.~S. Holling, S.~R. Carpenter, and A.~Kinzig, ``Resilience,
  adaptability and transformability in social--ecological systems,''
  \emph{Ecology and society}, vol.~9, no.~2, 2004.

\bibitem{steg2009encouraging}
L.~Steg and C.~Vlek, ``Encouraging pro-environmental behaviour: An integrative
  review and research agenda,'' \emph{Journal of environmental psychology},
  vol.~29, no.~3, pp. 309--317, 2009.

\bibitem{burke2006contemporary}
P.~J. Burke, \emph{Contemporary social psychological theories}.\hskip 1em plus
  0.5em minus 0.4em\relax Stanford University Press, 2006.

\bibitem{allison2001livelihoods}
E.~H. Allison and F.~Ellis, ``The livelihoods approach and management of
  small-scale fisheries,'' \emph{Marine policy}, vol.~25, no.~5, pp. 377--388,
  2001.

\bibitem{uribe2012community}
A.~L.~M. Uribe, D.~M. Winham, and C.~M. Wharton, ``Community supported
  agriculture membership in arizona. an exploratory study of food and
  sustainability behaviours,'' \emph{Appetite}, vol.~59, no.~2, pp. 431--436,
  2012.

\bibitem{bernstein09}
D.~S. Bernstein, \emph{Matrix Mathematics: Theory, Facts, and Formulas Ed.
  2}.\hskip 1em plus 0.5em minus 0.4em\relax Princeton University Press, 2009.

\bibitem{berman94}
A.~Berman and R.~J. Plemmons, \emph{Nonnegative matrices in the mathematical
  sciences}.\hskip 1em plus 0.5em minus 0.4em\relax Siam, 1994, vol.~9.

\end{thebibliography}
\end{document}